\documentclass[a4paper]{article}
\usepackage[english]{babel}
\usepackage[cp1251]{inputenc}
\usepackage[colorlinks,linkcolor=blue,citecolor=blue]{hyperref}
\usepackage{latexsym,amsfonts,amssymb,amsthm,amsmath,amsbsy,graphicx}
\begin{document}
\thispagestyle{myheadings} \markboth{Matrix factorization
identity}{Theory of Stochastic Processes, Vol. 11 (27), no. 1-2,
2005, pp.40-47}
\bigskip
 {\noindent \Large\bf\sc    Matrix factorization identity for almost semi-continuous\\ processes on a Markov
chain}\footnotemark[1] \footnotetext{This is an electronic reprint
of the original article published in Theory of Stochastic Processes,
Vol. 11 (27), no. 1-2, 2005. This reprint differs from the original
in pagination and typographic detail.}

\bigskip
\bigskip
 { \bf D.V. Gusak,\footnote{Institute of Mathematics Ukrainian National Academy of Science, 3 Tereshenkivska str.,
252601 Kiev, Ukraine. \phantom{\quad \quad }
 \href{mailto:random@imath.kiev.ua}{random@imath.kiev.ua}
 } E.V.~Karnaukh
 \footnote{Department of Probability and Mathematical Statistics, Kyiv National
University, 64~Vladimirs\-kaya str., 252017 Kyiv, Ukraine.
\phantom{\quad \quad }
 \href{mailto:kveugene@mail.ru}{kveugene@mail.ru}
 }}\hskip 6 cm UDC 519.21
\begin{center}
\bigskip
\begin{quotation}
\noindent  {\small In this article almost semi-continuous processes
with stationary independent increments on a finite irreducible
Markov chain are considered. For these processes the components of
matrix factorization identity are concretely defined. On the basis
of this concrete definition the relations for the distributions of
extrema and distributions of their complements for the almost upper
semi-continuous processes are established.}
\end{quotation}
\end{center}
\bigskip

\par The processes with stationary independent increments on a Markov
 chain are considered
in~[1]~-~[4]. These processes also were considered in~[7], where
they were called as in~[1,2] the risk processes in a Markovian
environment or Markov additive processes.
\par We'll consider the two-dimensional Markov process:
$$Z(t)=\{\xi(t),x(t)\}\quad\left( t\ge 0,\xi(0)=0\right),$$
where $x(t)$ is a finite ergodic Markov chain with state space
${E}'=\{1 \ldots m\}$ and transition matrix $ \mathbf{P}(t)=e^{t
\mathbf{Q}},$ $t\ge0,$ $
\mathbf{Q}=\mathbf{N}(\mathbf{P}-\mathbf{I}), $
where $\mathbf{N}=||\delta_{kr}\,\nu_k||^{m}_{k,r=1}$, $\nu_k$ are
parameters of the exponentially distributed random variables $
\zeta_k$ (the sojourn time of $x(t)$ in state $k$),
$\mathbf{P}=\|p_{kr}\|$ is the transition matrix of the embedded
Markov chain.
\par Let $\sigma_n=\sum _{k\leq n}\zeta _k$, $y_n=x(\sigma_n)$.
We suppose that under conditions
$x(\sigma_n-0)=k,$ $x(t)=r,$ $t\in[\sigma_n,\sigma_{n+1})$ $\xi (t)$
is determined by the processes
with stationary independent increments $\xi_r(t)\left(\xi_r(0)=0\right)$
 and by the independent
jumps $\chi_{kr}\,(k\neq r)$ at the moments $\sigma_n$. $\xi_r(t)$ have
the cumulant functions:
$$\psi_r(\alpha)=\imath\alpha a_r-\frac{1}{2}b_r^2\alpha^2+
\int_{ -\infty}^{\infty}[e^{\imath\alpha x}-1-\imath\alpha
x\delta(|x|\leq 1)]
 \Pi_r(dx),$$ $|a_r|<\infty, b_r^2<\infty,
\Pi_r(\cdot)$ are spectral measures of $\xi _r(t)$. And we denote
$\mathbf{F}(x)=\|\mathrm{P}\{\chi_{kr}<x\, ;y_1=r/y_0=k\}\|.$

\par The evolution of the process $Z(t)$ is determined by the matrix
characteristic function(ch.f.):
$$\boldsymbol{\Phi}_t(\alpha)=
\|\mathrm{E}[e^{\imath \alpha (\xi(t+u)-\xi(u))},x(t+u)= r/x(u)=k]
\|\qquad u\geq 0,$$ this ch.f. we can represent in the next form
$$\boldsymbol{\Phi}_t(\alpha)=\mathbf{E}e^{\imath \alpha \xi(t)}=e^{ t
\boldsymbol{\Psi}(\alpha)},\qquad \boldsymbol{\Psi}(0)=\mathbf{Q},$$
$$ \boldsymbol{\Psi}( \alpha)=\|\psi _k(\alpha)\delta _{kr}\|+
\mathbf{N}\left[\int_{-\infty}^{\infty}e^{\imath \alpha
x}d\mathbf{F} (x)-\mathbf{I}\right],
\boldsymbol{\Psi}(0)=\mathbf{Q}.
$$
\par Let $\theta _s$ denote an exponentially distributed random
variable with the parameter $s>0$ $(\mathrm{P}\{\theta _s>t\}=e^{-s
t}, t\ge 0 )$. We assume that $\theta _s$ is independent on $Z(t)$,
then
$$
\boldsymbol{\Phi}(s,\alpha )=\mathbf{E}e^{\imath \alpha \xi(\theta
_s)}= s\int_0^\infty e^{-st}\boldsymbol{\Phi}_t(\alpha
)dt=s\left(s\mathbf{I}-\boldsymbol{\Psi}(\alpha )\right)^{-1}.\eqno
(1)
$$
$$\mathbf{P}_s=s\int_0^{\infty}  e^{-st}\mathbf{P}(t)dt=
s\left(s\mathbf{I}-\mathbf{Q}\right)^{-1}.$$
\par Let us denote the next functionals for $\xi (t)$:
$$
\xi ^+(t)=\sup\limits_{0\leq u\leq t}\xi (u),\quad  \overline{\xi
}(t)=\xi (t)-\xi ^+(t),\quad \xi ^-(t)=\inf\limits_{0\leq u\leq
t}\xi (u), \quad \check{\xi} (t)=\xi (t)-\xi ^-(t),$$
$$ \tau^+(x)=\inf\{t:\xi (t)>x\},\;  \gamma
^+(x)=\xi (\tau ^+(x))-x,$$
$$
 \gamma _+(x)=x-\xi (\tau ^+(x)-0),\gamma _x^+ =
  \gamma ^+(x)+\gamma_+(x),\quad x>0,
 $$
 $$
 \tau^-(x)= \inf\{t:\xi (t)<x\},\quad x<0.
$$
and the distributions
$$\overline{\mathbf{P}}_+(s,x)=\mathbf{P}\left\{\xi ^+(\theta_s)>x\right\},x>0,
\mathbf{p}_+(s)=\mathbf{P}\left\{\xi ^+(\theta_s)=0\right\},
\mathbf{q}_+(s)=\mathbf{P}_s-\mathbf{p}_+(s);$$
$$\overline{\mathbf{P}}^+(s,x)=\mathbf{P}\left\{\check{\xi}(\theta_s)>x\right\},x>0,
\check{\mathbf{p}}_+(s)=\mathbf{P}\left\{\check{\xi}
^+(\theta_s)=0\right\},
\check{\mathbf{q}}_+(s)=\mathbf{P}_s-\check{\mathbf{p}}_+(s);$$
$$\mathbf{P}^-(s,x)=\mathbf{P}\left\{\overline{\xi}(\theta_s)<x\right\},
\mathbf{P}_-(s,x)=\mathbf{P}\left\{{\xi}^-(\theta_s)<x\right\},\;x<0.$$
{\bf Lemma 1}. [3] {\it For the two-dimensional Markov process
$Z(t)=\{\xi(t),x(t)\}$ the basic factorization identity is valid
$$
\boldsymbol{\Phi}\left(s,\alpha \right)=\mathbf{E}e^{\imath\alpha
\xi (\theta _s)}= \begin{cases}
\boldsymbol{\Phi}_+(s,\alpha )\mathbf{P}_s^{-1}\boldsymbol{\Phi }^-(s,\alpha ),\\
\boldsymbol{\Phi}_-(s,\alpha )\mathbf{P}_s^{-1}\boldsymbol{\Phi
}^+(s,\alpha ),
\end{cases}
\quad \text{where}\eqno(2) $$
$$
\boldsymbol{\Phi}_+(s,\alpha ) =\mathbf{E}e^{\imath\alpha
\xi^+(\theta_s)},\;
\boldsymbol{\Phi}^-(s,\alpha)=\mathbf{E}e^{\imath\alpha
\overline{\xi}(\theta_s)},\; \boldsymbol{\Phi}^+(s,\alpha
)=\mathbf{E}e^{\imath\alpha \check{\xi}(\theta_s)},$$
$$
\boldsymbol{\Phi}_-(s,\alpha )=\mathbf{E}e^{\imath\alpha
\xi^-(\theta_s)}.
$$
}
In papers~[3, 4] the concrete definition of components of (2) for
the semi-continuous processes was obtained. We want to derive
analogical concrete definition for the almost semi-continuous
processes, that were investigated for the scalar case
$\left(m=1\right)$ in paper~[5]. Before this let us consider some
auxiliary definitions and statements.
\par Let $\mathrm{B}_m(\alpha )$ is the
Banach algebra of matrices with dimension $m\times m$, the
elements of these matrices are the Fourier-Stieltjes transforms
 of the functions $f_{kr}(x)$, $k,r=1..m$, with bounded variation.
\par If $\Phi =\Phi ( \alpha )\in \mathrm{B}_m(\alpha ),
 \Phi (\alpha)=\|\int_{-\infty}^{\infty}e^{\imath \alpha x}df_{kr}(x)\|$,
  then we determine projection operations by the next relations
$$
\left[\Phi (\alpha )\right]_{\pm}=\|\pm \int_{\pm
0}^{\pm\infty}e^{\imath \alpha x}df_{kr}(x)\|,$$
$$ \left[\Phi (\alpha
)\right]_{\pm}^0=\|c_{kr}^{\pm}\pm \int_{\pm 0}^{\pm
\infty}e^{\imath \alpha x}df_{kr}(x)\|.
$$
\par The relation between the distribution
of $\xi^+ (\theta _s)$ and generating function of $\tau^+(z)$ is
valid:
$$ \overline{\mathbf{P}}_+(s,z)=\mathbf{E}\left [e^{-s
\tau ^+(z)},\tau^+(z)<\infty\right ]\mathbf{P}_s.\eqno(3)
$$%
In the sequel we will denote $\mathbf{E}\left [e^{-s \tau
^{\pm}(z)},\tau^{\pm}(z)<\infty\right ]=\mathbf{T}^{\pm}(s,x)$,
taking into account that the generating function of $\tau ^{\pm}(z)$
is considered on the chain $x(\tau ^{\pm}(z))$. Note that the state
space of the chain $x\left(\tau ^{\pm}(z)\right)$ can be narrow due
to nonattainability the level $z>0\,(z<0)$ by the process $\xi(t)$ .
Therefore we impose the next conditions
$$\{k:\mathrm{P}\{x\left(\tau^{\pm}(z)\right)=k\}>0\}={E}',
\{k:\mathrm{P} \{x\left(\tau^{\pm}(z)\right)=k\}=0\}=\varnothing.$$
$$
\left |\left |\,\mathrm{E}\left [\left |\xi (t)\right
|,x(t)=r/x(0)=k\right ]\,\right |\right |<\infty.\eqno(4)
$$ We'll use
the next notation
$$
d\mathbf{K}_0(z)=\mathbf{N}d\mathbf{F}(z)+\boldsymbol{\Pi}(dz),\,
\overline{\mathbf{K}}_0(x)=\int_x^{\infty}d\mathbf{K}_0(z),$$
$$ \mathbf{w}_+(\alpha
,u,v,\mu)=\int_0^{\infty}e^{\imath\alpha
x}\mathbf{W}(x,u,v,\mu)dx=\int_0^{\infty}e^{\imath\alpha
x}\int_x^{\infty}e^{(u-v)x-(u+\mu)z}d\mathbf{K}_0(z)dx ,
$$
\begin{multline*}\mathbf{v}^+(s,\alpha ,u,v,\mu )=\int_0^{\infty}e^{\imath\alpha x}
\mathbf{V}_+ (s,x,u,v,\mu )dx\\=\int_0^{\infty}e^{\imath\alpha
x}\mathbf{E}\left[e^{-s\tau ^+(x)-u\gamma ^+(x)-v\gamma
_+(x)-\mu\gamma _x^+},\tau ^+(x) <\infty\right]dx,
\end{multline*}
$${C}_*(s)=\begin{cases}
s^{-1}\left(\frac{1}{2}\frac{\partial}{\partial
x}\mathbf{P}^-(s,x)|_{x=0}\mathbf{B}^2+
\left[\boldsymbol{\Phi}^-(s,\alpha)\right]_+^0\mathbf{A}_0^+\right),&
\mathbf{B} \geq\mathbf{O};\\
s^{-1} \left[\boldsymbol{\Phi}^-(s,\alpha)\right]_+^0\mathbf{A}^+,
&\mathbf{B}\equiv \mathbf{O}.
\end{cases}$$
$$\mathbf{A}_0^+=\|\delta _{kr}\delta (b_k=0,a_k>0)a_k\|,\,
 \mathbf{A}^+=\|\delta _{kr}\delta(a_k>0)a_ k\|, \mathbf{B}=
 \|\delta _{kr}b^2_r\|.$$
{\bf Lemma 2}. [3] {\it If condition~{\rm{(3)}} is satisfied, then
$$
{\mathbf{v}^+(s,\alpha ,u,v,\mu
)=\boldsymbol{\Phi}_+(s,\alpha)\mathbf{P}_s^{-1}\left({C}_*(s)+s^
{-1}\left[\boldsymbol{\Phi}^-(s,\alpha)\mathbf{w}_+(\alpha ,u,v,\mu
)\right]_+^0\right)}\eqno(5)
$$
}
If we denote
$${\mathbf{K}(s,x)=\int^0_{-\infty}d\mathbf{P}^-(s,y)
\overline{\mathbf{K}}_0(x-y),\quad \mathbf{k}(s,\alpha
)=\int_0^{\infty}e^{i\alpha x}\mathbf{K}(s,x)dx ,}
$$
then from lemma 2 the next statement follows.

{\bf Theorem 1}.{\it If condition~{\rm{(4)}} is satisfied, then
$${\boldsymbol{\Phi } _+(s,\alpha
)=\left(\mathbf{I}-\imath\alpha\left({C}_*(s)+s^{-1}\mathbf{k
}(s,\alpha )\right)\right)^{-1}\mathbf{P}_s.}\eqno (6)
$$
}
\begin{proof} Let's substitute $u=v=\mu =0$ in formula~(5), then
$${\mathbf{ v}^+(s,\alpha ,0,0,0 )=\boldsymbol{\Phi}
_+(s,\alpha)\mathbf{P}_s^{-1}\left({
C}_*(s)+s^{-1}\left[\boldsymbol{\Phi
}^-(s,\alpha)\mathbf{w}_+(\alpha ,0,0 ,0
)\right]_+^0\right),}\eqno(7)
$$
where
$$\left[\boldsymbol{\Phi }^-(s,\alpha )\mathbf{w}_+(\alpha,0,0,0
)\right]_+^0=\left[\int_{-\infty}^0 e^{\imath\alpha
x}d\mathbf{P}^-(s,x) \int_0^{\infty}e^{\imath\alpha
x}\overline{\mathbf{K}}_0(x)dx\right]_+^0$$
$$
=\int_0^{\infty}e^{i\alpha
x}\int_{-\infty}^{0}d\mathbf{P}^-(s,y)\overline{\mathbf{K}}_0(x-y)dx=
\mathbf{k}(s,\alpha ).
$$
Using formula~(3) and definition of $\mathbf{v}^+$, we have
\begin{multline*}
\mathbf{ v}^+(s,\alpha ,0,0,0)=\int_0^{\infty}e^{\imath \alpha x}
\mathbf{E}\left [e^{ -s\tau ^+(x)},\tau^+(x)<\infty\right ]dx=
\int_0^{\infty}e^{\imath \alpha x}\mathbf{ P}\{\xi
^+(\theta_s)>x\}dx\,\mathbf{P}_s^{-1}\\
=\frac{1}{\imath \alpha }\left(\boldsymbol{\Phi }_+(s,\alpha
)-\mathbf{P}_s\right)\mathbf{P}_s^{-1}.
\end{multline*}
Let's substitute received relations in formula~(7):
$$\frac{1}{\imath\alpha }\left(\boldsymbol{\Phi }_+(s,\alpha )
-\mathbf{P}_s\right)\mathbf{P}_s^{-1} =\boldsymbol{\Phi}_+(s,\alpha
)\mathbf{P}_s^{-1} \left({C}_*(s)+s^{-1}\mathbf{k}( s,\alpha
)\right),$$
whence we receive{{~\rm{(6)}}}. In paper~[3] the representation of
$\boldsymbol{\Phi} _+(s,\alpha )$ in terms of the joint generating
function of $\{\tau ^+(0),\gamma ^+(0)\}$ was obtained.
Formula{{~\rm{(6)}}} is the another way of representation of
$\boldsymbol{\Phi} _+(s,\alpha )$ in terms of the integral transform
of the convolution $\mathbf{P}^-(s,x)$ with
$\overline{\mathbf{K}}_0(x).$ In some partial cases the distribution
of $\mathbf{P}^-(s,x)$ has the exponential form. Then the integral
transform of the convolution $\mathbf{k}(s,\alpha)$ has more simple
form.
\end{proof}
Let us consider the analogy of the almost upper semi-continuous
scalar process analyzed in paper~[5].
The process $Z(t)=\{\xi(t),x(t)\}$ is the almost upper
semi-continuous process on a Markov chain, if
$\mathbf{B}\equiv\mathbf{O},\,$ $\mathbf{A}<0$, and
${\int_{-\infty}^{\infty}\boldsymbol{\Pi}(dx)=\boldsymbol{\Lambda
}<\infty}$, $\boldsymbol{\Pi}(dx)=\boldsymbol{\Lambda
}d\mathbf{F}_0(x)$, where $\mathbf{F}_0(x)=\|\delta
_{kr}F_0^k(x)\|,\,F_0^k(x)$ are the distribution functions of jumps
of $\xi (t)$, when $x(t)=k$. $\boldsymbol{\Lambda}=\|\delta
_{kr}\lambda _k\|$, $\lambda _k$ are parameters of the exponentially
distributed random variables $\zeta'_k$ (the time between two
consecutive jumps of $\xi (t) $ if $x(t)=k$). We also assume that
$d\mathbf{K}_0(z)= \boldsymbol{\Lambda
}\overline{\mathbf{F}}_0(0)\mathbf{C}e^{-\mathbf{C}z}dz,\, z>0$,
where $\mathbf{C}=\|\delta _{kr}c_{k}\| (c_k>0)$ (the positive jumps
of $\xi (t)$ have exponential distributions with the parameters
$c_k$ if $x(t)=k$). Under these conditions the cumulant of $Z(t)$
has the next form
$$\boldsymbol{\Psi }(\alpha )=\imath\alpha \mathbf{A}+\boldsymbol{\Lambda }
\overline{\mathbf{F}}_0(0)\mathbf{C}\left(\mathbf{C}- \imath\alpha
\mathbf{I}\right)^{-1} +\int^0_{-\infty}e^{\imath\alpha
x}d\mathbf{K}_0(x) -\boldsymbol{\Lambda}-\mathbf{N}.$$
For the almost upper semi-continuous process $Z(t)$ we have the next
concrete definition of components of the first part of{{~\rm{(2)}}}.

{\bf Theorem 2}.{\it For the almost upper semi-continuous processes
the distribution of $\xi ^+(\theta _s)$ is determined by the next
relations
$$
{\boldsymbol{\Phi }_+(s,\alpha)=\left(\mathbf{C}-\imath\alpha
\mathbf{I}\right)\left(\mathbf{p}_+(s)\mathbf{P}_s^{-1}\mathbf{C}-\imath\alpha
\mathbf{I}\right)^{ -1}\mathbf{p}_+(s),}\eqno(8)
$$
$${\mathbf{p}_+(s)=s\left(s\mathbf{I}+\mathbf{E} \,e^{\overline{\xi}
(\theta _s)\mathbf{C}}
\boldsymbol{\Lambda}\overline{\mathbf{F}}_{0}(0)\right)^{-1}
\mathbf{P}_s;}\eqno(9)
$$
$${\overline{\mathbf{P}}_{+}(s,x)=\mathbf{q}_{+}(s)e^{-\mathbf{P}_{s}^{-1}
\mathbf{C}\mathbf{p}_{+}(s)x},\quad x>0, }\eqno(10)
$$
for $\overline{\xi }(\theta _s)=\xi (\theta_s)-\xi ^+(\theta_s)$ we
have
$${\boldsymbol{\Phi }^{-}(s,\alpha )=
\mathbf{P}_s\mathbf{p}_+^{-1}(s)\left [\boldsymbol{\Phi}(s,\alpha
)\right ]_- -\mathbf{q}_{+}(s)\mathbf{p}_+^{-1}(s)\mathbf{C}\left
[\left(\mathbf{C}-\imath\alpha\mathbf{I}\right)^{-1}
\boldsymbol{\Phi}(s,\alpha)\right ]_{-}, }\eqno(11)$$
$${\mathbf{P}^{-}(s,x)=\mathbf{P}_s\mathbf{p}_+^{-1}(s)\mathbf{P}(s,x)
-\mathbf{q}_{+}(s)\mathbf{p}_+^{-1}(s)\mathbf{C}\int_0^{\infty}
e^{-\mathbf{C}y}\mathbf{P}(s,x-y)dy\,,\quad x<0 .}\eqno(12)
$$
}

\begin{proof} Substituting $d\mathbf{K}_0(z)= \boldsymbol{\Lambda
}\overline{\mathbf{F}}_0(0)\mathbf{C}e^{-\mathbf{C}z}dz,\, z>0$ in
formula~(6), we obtain
$${ \boldsymbol{\Phi }_+(s,\alpha )= \left(\mathbf{I}-\imath\alpha\left({C}_*
({s})+s^{-1}\int^0_{-\infty}d\mathbf{P}^-(s,y)e^{\mathbf{C}y}
\boldsymbol{\Lambda}\overline{\mathbf{F}}_0(0)\left(\mathbf{C}-\imath\alpha
\mathbf{I}\right)^{-1}\right)\right)^{-1 }\mathbf{P}_s.}\eqno(13)
$$
Under conditions of the theorem: ${C}_*(s)\equiv\mathbf{O}$, then
$${\boldsymbol{\Phi }_+(s,\alpha )=\left(\mathbf{I}-\imath\alpha s^{-1}
\int^0_{-\infty
}d\mathbf{P}^-(s,y)e^{\mathbf{C}y}\boldsymbol{\Lambda}
\overline{\mathbf{F}}_0(0)
\left(\mathbf{C}-\imath\alpha\mathbf{I}\right)^{-1}\right)^{-1}
\mathbf{P}_s.}\eqno(14)
$$
Hence
$${\mathbf{p}_+(s)=\lim_{\imath\alpha\rightarrow\infty}\boldsymbol{\Phi
}_+(s,\alpha
)=\left(\mathbf{I}+s^{-1}\int^0_{-\infty}d\mathbf{P}^-(s,y)e^{\mathbf{C}y}
\boldsymbol{\Lambda}\overline{\mathbf{F}}_0(0)\right)^{-1}\mathbf{P}_s.}\eqno(15)
$$
Substituting (15) in (14) and taking into account that
$\int^0_{-\infty}d\mathbf{P}^-(s,y)e^{\mathbf{C}y}=
\mathbf{E}\,e^{\overline{\xi}(\theta _s)\mathbf{C}}$ we
obtain{{~\rm{(8)}}}.
 To prove formula~(10) let us invert~(5) with respect to $\alpha$:
$${s\mathbf{V}_+(s,x,u,v,\mu)=\int_{0}^{x}d\mathbf{P}_+(s,z)
\mathbf{P}_s^{-1}\int_{-\infty}^0
d\mathbf{P}^-(s,y)\mathbf{W}(x-y-z,u,v,\mu).}\eqno (16)
$$
Letting  $u=v=\mu =0$ we have
$$s\mathbf{V}_+(s,x,0,0,0)=\int_{0}^{x}d\mathbf{P}_+(s,z)
\mathbf{P}_s^{-1}\int _{-\infty}^0 d\mathbf{P}^-(s,y)
\overline{\mathbf{K}}_0(x-y-z).$$
Taking into account{{~\rm{(3)}}} and the condition of the almost
semi-continuity we obtain the equation
$${s\mathbf{T}^+(s,x)=-\int_{0}^{x}d\mathbf{T}^+(s,z)\int_{-\infty} ^0
d\mathbf{P}^-(s,y)\boldsymbol{\Lambda}
\overline{\mathbf{F}}_0(0)e^{\mathbf{C}(y+z)}e^
{-\mathbf{C}x}.}\eqno(17)
$$
Let's differentiate with respect to $x>0$ then we have the next
equation for $\mathbf{T}^+(s,x)$:
$$
\displaystyle\frac{\partial}{\partial
x}\mathbf{T}^+(s,x)=-\mathbf{T}^+(s,x)\mathbf{C}\mathbf{p}_+(s)
\mathbf{P}_s^{-1}\,,\quad x>0,\eqno (18)$$
$$
\displaystyle\mathbf{T}(s,0)=\mathbf{q}_+(s)\mathbf{P}_s^{-1}.
$$
The solution of equation{{~\rm{(18)}}} is expressed by the next
formula
$${ \mathbf{T}^+(s,x)=\mathbf{T}(s,0)e^{-\mathbf{C}\mathbf{p}_+
(s)\mathbf{P}_s^{-1}x}.}\eqno(19)
$$
Taking into account formulas{{~\rm{(19)}}} and{{~\rm{(3)}}},
 we have
\begin{multline*}
\overline{\mathbf{P}}_{+}(s,x)=\mathbf{T}^{+}(s,x)\mathbf{P}_s\\
=\mathbf{q}_{+}(s)\mathbf{P}_s^{-1}e^{-\mathbf{C}\mathbf{p}_{+}(s)
\mathbf{P}_{s}^{-1}x}\mathbf{P}_s\\
=\mathbf{q}_{+}(s)e^{-\mathbf{P}_{s}^{-1}\mathbf{C}\mathbf{p}_{+}(s)x}.
\end{multline*}
The representations of $\boldsymbol{\Phi }^{-}(s,\alpha )$
in{{~\rm{(11)}}} we can obtain from the first part of factorization
identity{{~\rm{(2)}}}, substituting instead of $\boldsymbol{\Phi
}_+(s,\alpha )$ its representation from formula{{~\rm{(8)}}}.
Formula{{~\rm{(12)}}} is received by the inversion of{{~\rm{(11)}}}
with respect to $\alpha$ and by the integration with respect to $z$.
\end{proof}
Further we consider the concrete definition of components of the
second part of{{~\rm{(2)}}}.

{\bf Theorem 3}. {\it For the almost upper semi-continuous processes
the distribution of $\check{\xi}(\theta_s)=\xi (\theta_s)-\xi
^-(\theta_s)$ is determined by the next relations:
$${\boldsymbol{\Phi
}^+(s,\alpha)=\check{\mathbf{p}}_+(s)\left(\mathbf{C}\mathbf{P}_s^{-1}
\check{\mathbf{p}}_+(s)-\imath\alpha \mathbf{I}\right)^{
-1}\left(\mathbf{C}-\imath\alpha \mathbf{I}\right),}\eqno(20)
$$
$$
{\check{\mathbf{p}}_+(s)=s\mathbf{P}_s\left(s\mathbf{I}+
\boldsymbol{\Lambda}\overline{\mathbf{F}}_{0}(0)\mathbf{E}
\,e^{\mathbf{C}{\xi}^- (\theta _s)} \right)^{-1},}\eqno(21)
$$
$${\overline{\mathbf{P}}^+(s,x)=e^{-\check{\mathbf{p}}_{+}(s)
\mathbf{C}\mathbf{P}_{s}^{-1}x} \check{\mathbf{q}}_{+}(s),\quad x>0;
}\eqno(22)
$$
for the minimum ${\xi}^-(\theta _s)$ we have:
$${\boldsymbol{\Phi }_-(s,\alpha)= \left [\boldsymbol{\Phi}(s,\alpha
)\right ]_-\check{\mathbf{p}}_+^{-1}(s)\mathbf{P}_s -\left
[\boldsymbol{\Phi}(s,\alpha)\left(\mathbf{C}-
\imath\alpha\mathbf{I}\right)^{-1}\right
]_{-}\mathbf{C}\check{\mathbf{p}}_+^{-1}(s)\check{\mathbf{q}}_{+}(s),
}\eqno(23)
$$
$${\mathbf{P}_{-}(s,x)=\mathbf{P}(s,x)\check{\mathbf{p}}_+^{-1}(s)
\mathbf{P}_s -\int_0^{\infty}
\mathbf{P}(s,x-y)e^{-\mathbf{C}y}dy\mathbf{C}\check{\mathbf{p}}_+^{-1}(s)
\check{\mathbf{q}}_{+}(s)\,,\quad x<0 .}\eqno(24)
$$
}
\begin{proof} Note that we'll consider such trajectories of the
process for which $\{\tau^-(x)<\infty\}$. Then the stochastic
relations for $\tau _ {kr}^-(x)$ $(k,r=\overline{1,m})$ if $x<0$,
where lower indices denote the initial state and the state of $x(t)$
at the moment of achievement the level $x$, correspondingly
$\left(x(0)= k,\,x\left(\tau^-(x)\right)=r\right)$, have the next
form

$${ \mathbf{\tau}_{kr}^-(x)\doteq \begin{cases}
x/{a_{k}}&\zeta _{k}'>x/{a_{k}},\quad\zeta _{k}>x/{a_{k}};\\
\zeta _k^{'}& \xi _k+a_k\zeta _k^{' }<x,\quad\zeta _k^{'}<\zeta _k;\\
\zeta _k & \chi _{kr}+a_k\zeta _k<x,\quad\zeta _k^{'}>\zeta _k;\\
\zeta _k^{'}+\tau _{kr}^+\left( x-a_k\zeta _k^{'}-\xi _k\right)&
\xi_k+a_k\zeta_k^{'}>x,\quad\zeta _k^{'}<\zeta_k;\\
\zeta _k+\tau
_{jr}^+\left(x-a_k\zeta_k-\chi_{kj}\right)&\chi_{kj}+a_k\zeta_k>x,
\quad\zeta_k^{'}>\zeta_k;
\end{cases}}\eqno(25)
$$
On the basis of{{~\rm{(25)}}}, we derive the equation
$$
\displaystyle T_{kr}^-(s,x)=\mathrm{E}\left [e^{-s\tau ^-(x)},\tau
^-(x)<\infty,\,x(\tau ^-(x))=r/x(0)=k\right ]$$
$$=\displaystyle
\delta_{kr}e^{-\left(s+\lambda _k+\nu _k\right) x/a_k}
+\frac{1}{a_k}\int_{-\infty}^{x}e^{-\left(\lambda _k
+\nu_k+s\right)\frac{x-y}{a_k}}\int_{-\infty}^{y}d{K}^{0}_{kr}(z)dy$$
$$
\displaystyle +\frac{1}{a_k}\sum_{j=1}^{m}\int_{-\infty}^x
e^{-\left(\lambda _k
+\nu_k+s\right)\frac{x-y}{a_k}}\int^{\infty}_{y}dK_{kj}^0(z)T_{jr}^-(s
,y-z)dy.\eqno(26)
$$
Let us differentiate left and right side of formula{{~\rm{(26)}}}
with respect to $x<0$:
$$ a_k\frac{\partial}{\partial
x}T_{kr}^-(s,x)=-\left(\lambda_k+\nu_k+s\right)T_{kr}^-(s,x)+$$
$$
+\int_{-\infty}^{x}dK^{0}_{kr}(z)+\sum_{j=1}^{m}\int^{\infty}_{x}
dK_{kj}^0(z)T_{jr}^-(s, x-z). \eqno(27)$$
Integro-differential equations (27) we can represent in the matrix
form:
$$ \displaystyle \mathbf{A}\frac{\partial}{\partial
x}\mathbf{T}^-(s,x)=-\left(s\mathbf{I}+\mathbf{N}+
\boldsymbol{\Lambda}\right)\mathbf{T}^-(s,x)+\int_{-\infty}^{x}
d\mathbf{K}_{0}(z)+$$
$$
+\int^{\infty}_{x}d\mathbf{K}_0(z)\mathbf{T}^-(s,x-z),\,x<0.\eqno(28)
$$
Analogically to{{~\rm{(3)}}} we have the formula of the relation
between the distribution of $\xi^- (\theta _s)$ and generating
function of $\tau^-(z)$:
$${\mathbf{T}^{-}(s,x)=\mathbf{I}-\overline{\mathbf{P}}_-(s,x)
\mathbf{P}^{-1}_s,\qquad
\overline{\mathbf{P}}_-(s,x)=\mathbf{P}\{\xi ^-(\theta
_s)>x\}.}\eqno(29)
$$
Substituting{{~\rm{(28)}}} in{{~\rm{(29)}}}, we obtain the next
equation
$${-\mathbf{A}\frac{\partial}{\partial
x}\overline{\mathbf{P}}_-(s,x)=\left(s\mathbf{I}+\mathbf{N}
+\boldsymbol{\Lambda}\right)\overline{\mathbf{P}}_-(s,x)-\int^{\infty}_{x}
d\mathbf{K}_0(z)\overline{\mathbf{P}}_-(s,x-z)-s\mathbf{I},\;x<0.}\eqno(30)
$$
Let's consider the Laplace transform of equation{{~\rm{(30)}}} with
respect to $x<0$:
$$\left(s\mathbf{I}-\boldsymbol{\Psi}(-\imath u)\right)
\boldsymbol{\Phi}_-(s,-\imath u)=s\mathbf{I}+$$
$$+\int_0^{\infty}
d\mathbf{K}_0(z)\int_{-\infty}^{0}\left(e^{ux}\overline{
\mathbf{P}}_{-}(s,x-z)-e^{(x+z)u}\overline{\mathbf{P}}_-(s,x)\right)dx.$$
Taking into account that $d\mathbf{K}_0(z)= \boldsymbol{\Lambda
}\overline{\mathbf{F}}_0(0)\mathbf{C}e^{-\mathbf{C}z}dz,\, z>0$, we
have:
$${\left(s\mathbf{I}-\boldsymbol{\Psi}(-\imath u)\right)
\boldsymbol{\Phi}_-(s,-\imath u)=s\mathbf{I}-u\boldsymbol{\Lambda
}\overline{\mathbf{F}}_0(0)\left(\mathbf{C}-u\mathbf{I}\right)^{-1}
\mathbf{E}\,e^{\mathbf{C}\xi ^{-}(\theta _s)}.}\eqno(31)
$$
Combining formula{{~\rm{(1)}}} with the second row of{{~\rm{(2)}}}
we obtain from{{~\rm{(31)}}} that
$${\boldsymbol{\Phi }^+(s,-\imath u )=\mathbf{P}_s\left(\mathbf{I}-u
s^{-1}\boldsymbol{\Lambda}\overline{\mathbf{F}}_0(0)\left(\mathbf{C}
-u\mathbf{I}\right)^{-1}\mathbf{E}\,e^{\mathbf{C}\xi ^-(\theta
_s)}\right)^{-1}.}\eqno(32)
$$
After the limit passage $u\rightarrow\infty$ in{{~\rm{(32)}}}
formula{{~\rm{(21)}}} follows. Substituting{{~\rm{(21)}}}
in{{~\rm{(32)}}} we obtain{{~\rm{(20)}}}. After inversion
of{{~\rm{(20)}}} with respect to $\alpha$ formula{{~\rm{(22)}}}
follows. Substituting{{~\rm{(20)}}} in the second row
of{{~\rm{(2)}}} we obtain {{~\rm{(23)}}}. Formula{{~\rm{(24)}}} is
received by the inversion of{{~\rm{(23)}}} with respect to
$\alpha$
and by the integration with respect to $z$.
\end{proof}
\noindent {\bf{Remark 1}}. Let ${Z}(t)=\left\{ {\xi}(t),{x}(t)\right\}$
 be the almost upper
semi-continuous process, then the process
${Z}_1(t)=\left\{ {-\xi}(t),{x}(t)\right\}$ is the
almost lower semi-continuous process. Taking into account
the next relations between the ch.f. of
extrema for $Z(t)$ and $Z_1(t)$:
$$\boldsymbol{\Phi }_{(1)}^{\pm}(s,\alpha )=
\boldsymbol{\Phi }_{\pm}(s,-\alpha ),
\quad \boldsymbol{\Phi }^{(1)}_{\pm}(s,\alpha )=\boldsymbol{\Phi
}^{\pm}(s,-\alpha ),$$ we can obtain the statements about the
concrete definition of components of{{~\rm{(2)}}} for the almost
lower semi-continuous processes.
\par\smallskip
\noindent {\bf{Remark 2}}. In theorems 2 and 3 we have considered
the case $\mathbf{A}<0$. Analogical results take place if we assume
that $\mathbf{A}=0$. In this case we should take into account that
$\mathbf{p}_-(s)=\mathbf{P}\left\{\xi ^-(\theta _s)=0\right\}\neq 0$
and $\check{\mathbf{p}}_-(s)=\mathbf{P}\left\{\overline{\xi}(\theta
_s)=0\right\}\neq 0$.

We'll need some further notation. Let $ \left(\pi _1,\ldots,\pi
_m\right)$ be the stationary distribution of $x(t)$,
$$
m_1^0=\sum_{k=1}^m \pi_k\sum_{r=1}^m m_{kr},\quad m_{kr}=\delta
_{kr}\left(a_k+\int_{R}x\Pi_k(dx)\right)+\int_{R}x\,\nu _k
dF_{kr}(x),
$$
$$
 \mathbf{K}(r)=\boldsymbol{\Psi }(-\imath r ).
$$
By results of~[6] it follows that
$$
  \lim_{t\rightarrow \infty}\mathbf{P}(t)=\lim_{s \rightarrow 0}
  \mathbf{P}_{s}=\lim_{s \rightarrow
0}s\left(s\mathbf{I}-\mathbf{Q}\right)^{-1}=\mathbf{P}_0,$$
$$
\mathbf{P}_0=\|p_{kr}^0\|,\quad p_{kr}^0=\pi_r>0,
$$
if $\left |{m_1^0}\right |>0$ then
$$\lim_{r\rightarrow 0}r\mathbf{K}^{-1}(r)=
\lim_{r\rightarrow 0}r\left(\mathbf{Q}+r\mathbf{M}_1\right)^{-1}=
\frac{1}{m_1^0}\mathbf{P}_0,\quad \mathbf{M}_1=\mathbf{A}+\int_{R} x
 d\mathbf{K}_0(x).$$

{\bf Theorem 4}. {\it If $m_1^0>0$ then for the almost upper
semi-continuous process $Z(t)$ the distribution of strictly negative
values of $\xi ^-=\inf\limits_{0\leq u\leq \infty}\xi(u)$ is
determined by the generating function
$${\mathbf{E}\left [e^{r\xi ^-},\xi
^-<0\right ]=\left
[r\mathbf{K}^{-1}(r)\left(\mathbf{C}-r\mathbf{I}\right)^{-1}\right
]_-\check{\mathbf{R}}_+,}\eqno(33)
$$
if $\mathbf{A}=0$ then
$${\mathbf{p}_- = \mathbf{P}\left\{\xi ^-=0\right\}=
\left(\boldsymbol{\Lambda
}-\mathbf{N}\left(\mathbf{f}(0)-\mathbf{I}\right)\right)
\check{\mathbf{R}}_+,}\eqno(34)
$$
where $\check{\mathbf{R}}_+=\lim\limits_{s\rightarrow
0}s\check{\mathbf{R}}^{-1}_+(s)=\lim\limits_{s\rightarrow
0}s\check{\mathbf{p}}^{-1}_{+}(s)\mathbf{P}_s$,
$\mathbf{f}(0)=\|\mathrm{P}\left\{\chi
_{kr}=0,y_1=r/y_0=k\right\}\|$. }

\begin{proof} From formula{{~\rm{(20)}}} and factorization
identity{{~\rm{(2)}}} it follows that
$$\boldsymbol{\Phi }_-(s,\alpha )=
\boldsymbol{\Phi }(s,\alpha )\left(\mathbf{C}-\imath\alpha
 \mathbf{I}\right)^{-1}
\left(\mathbf{C}\mathbf{P}_s^{-1}\check{\mathbf{p}}_+(s)-
\imath\alpha \mathbf{I}\right)
\check{\mathbf{p}}^{-1}_{+}(s)\mathbf{P}_s.$$ or
$${\mathbf{E}e^{r\xi ^-(\theta
_s)}=s\left(s\mathbf{I}-\mathbf{K}(r)\right)^{-1}\left(\mathbf{C}-
r\mathbf{I}\right)^{-1}
\left(\mathbf{C}-r\check{\mathbf{R}}^{-1}_+(s)\right).}\eqno(35)
$$
From relation{{~\rm{(35)}}} it follows that
$$
\lim_{r\rightarrow 0}\lim_{s\rightarrow 0}\mathbf{E}e^{r\xi
^-(\theta_s)}=\frac{1}{m_1^0}\mathbf{P}_0\mathbf{C}^{-1}
\check{\mathbf{R}}_+.
$$
Hence the condition $m_1^0>0$ provides the existence of
$\lim\limits_{s\rightarrow
0}s\check{\mathbf{R}}^{-1}_+(s)=\;\check{\mathbf{R}}_+$. Then
from{{~\rm{(35)}}} after the limit passage $\left(s\rightarrow
0\right)$ formula{{~\rm{(33)}}} follows. If $\mathbf{A}=0$ let us
consider the probability $P^{0}_{kr}(t)=\mathrm{P}\left\{\xi
(t)=0,\,x(t)=r/x(0)=k\right\}$ which satisfies the next equation
$${P^{0}_{kr}(t)=\delta _{kr}e^{-\left(\nu _k+\lambda _k\right)t}
+\int_0^t \nu _k e^{-\left(\nu _k+\lambda
_k\right)y}\sum_{j=1}^{m}\mathrm{P}\left\{\chi_{kj}=0,\,y_1=j/y_0=k
\right\}P^{0}_{jr}(t-y)dy.}\eqno(36)
$$
Applying Laplace-Karson transform to equation{{~\rm{(36)}}} we
obtain the next equation
$${\widetilde{\mathbf{P}}^{0}(s)=\left(s\mathbf{I}+\mathbf{N}
+\boldsymbol{\Lambda
}\right)^{-1}\left(s\mathbf{I}+\mathbf{N}\mathbf{f}(0)\right)
\widetilde{\mathbf{P}}^{0}(s)}
$$
or
$${\widetilde{\mathbf{P}}^{0}(s)=s\left(s\mathbf{I}
+\boldsymbol{\Lambda
}-\mathbf{N}\left(\mathbf{f}(0)-\mathbf{I}\right)\right)^{-1}.}\eqno(37)
$$
Applying to the first part of{{~\rm{(2)}}} operation $\left
[\,\right ]^0$ and taking into account formula{{~\rm{(37)}}} we
obtain relation
$$\mathbf{p}_-(s)\mathbf{P}_s^{-1}\check{\mathbf{p}}_+(s)=\mathbf{p}_-(s)
\check{\mathbf{R}}_+(s)=s\left(s\mathbf{I}+\boldsymbol{\Lambda
}-\mathbf{N}\left(\mathbf{f}(0)-\mathbf{I}\right)\right)^{-1}$$
or
$${\mathbf{p}_-(s)=\left(s\mathbf{I}+\boldsymbol{\Lambda
}-\mathbf{N}\left(\mathbf{f}(0)-\mathbf{I}\right)\right)^{-1}\,s
\check{\mathbf{R}}^{-1}_+(s).}\eqno(38)
$$
After the limit passage $\left(s\rightarrow 0\right)$
from{{~\rm{(38)}}} formula{{~\rm{(34)}}} follows.
\end{proof}

\begin{thebibliography}{1}
%
\bibitem{Cinlar} E.~\c{C}inlar, "Markov additive processes",
Z. Wahrscheinlichkeitstheor. und verw. Geb., No.~24,~85-121 (1972).
%
\bibitem{Miller} H.D.~Miller, "A matrix factorization problem in
the theory of random variables defined on a finite Markov chain.
Absorption probabilities for sums of random variables defined on a
finite Markov chain", Proc. Cambrige Philos. Soc., {\bf 58}, No. 2,
 pp. 268-298 (1962).
%
\bibitem{Gusak2}D.V. Husak,
\emph{Boundary-Value Problems for Processes with Independent
Increments on Finite Markov Chains and for Semi-Markov Processes}
[in Ukrainian], Institute of Mathematics, Ukrainian Academy of
Sciences, Kyiv (1998).
%
\bibitem{Gusak3}  D.V.~Gusak, "A factorization identity for
semicontinuous processes defined on a Markov chain",  Theor.
Probability and Math. Statist., No. 64 , pp. 37-50 (2002).
%
\bibitem{Gusak4} D.V.~Gusak, "The connection of distributions of
extrema for Levy processes with its ladder points", Theory of
Stochastic
 processes, {\bf 10(26)}, 3-4, 35-42  (2004).
%
\bibitem{Korolyuk} V.S.~Korolyuk, A.F.~Turbin, \emph{Semi-Markov processes
and their applications}[in Russian], N.Dumka, Kiev (1976).
%
\bibitem{Asmussen} S.~Asmussen, \emph{Ruin Probabilities},  Word Scientist,
Singapore(2000).
\end{thebibliography}
\end{document}